\documentclass[11pt]{article}
\usepackage{amsmath}
\usepackage{amssymb}
\usepackage{amsthm}
\usepackage{epsfig}

\newtheorem{theorem}{Theorem}
\newtheorem{lemma}{Lemma}

\newtheorem{proposition}{Proposition}

\theoremstyle{remark}
\newtheorem{remark}{Remark}

\renewcommand{\Pr}{\mathbb P}
\renewcommand{\limsup}{\varlimsup}
\renewcommand{\liminf}{\varliminf}
\newcommand{\expect}{\mathbb E}
\newcommand{\e}{\varepsilon}
\newcommand{\ind}{1}
\renewcommand{\epsilon}{\varepsilon}
\newcommand {\be}[1]{\begin{equation}\label{#1}}
\newcommand {\ee}{\end{equation}}

\begin{document}

\title{A Transposition Rule Analysis Based on a Particle Process}

\author{David Gamarnik\\
\small IBM T.J. Watson Research Center\\
\small Yorktown Heights, NY 10598\\ \small gamarnik@watson.ibm.com \and Petar Mom\v cilovi\'c\thanks{The work was performed while with IBM T.J. Watson Research Center, Yorktown Heights, NY 10598.} \\
\small EECS Department \\
\small University of Michigan \\
\small Ann Arbor, MI 48109 \\
\small petar@eecs.umich.edu}

\maketitle

\begin{abstract}
A linear list is a collection of items that can be accessed sequentially. The cost of a request is the number of items that need to be examined before the desired item is located, i.e., the distance of the requested item from the beginning of the list. The transposition rule is one of the algorithms designed to reduce the search cost by organizing the list. In particular, upon a request for a given item, the item is transposed with the preceding one. We develop a new approach for analyzing the algorithm. The approach is based on a coupling with a certain constrained asymmetric exclusion process. This allows us to establish an asymptotic optimality of the rule for two families of request distributions.
\end{abstract}


\section{Introduction}

The linear list, a collection of items that can be accessed sequentially, is one of basic data structures known in computer science. A primary operation defined on the list is search. A requested item is found in the list by sequentially examining items from the beginning of the list one by one. The cost of search is defined to be the distance of the requested item from the beginning of the list, i.e., the number of items that need to be examined in order to locate the desired item. Intuitively, one would like to place frequently requested items at the front of the list so that to minimize the number of search steps. If the request sequence were known a priori, one could place items in an order that minimizes the search cost. Yet often properties of the request sequence are either not known in advance or time dependent. Hence, it is desirable to employ an algorithm that organizes the list based on past requests. The two best known self-organizing algorithms are the move-to-front rule and transposition rule~\cite[Section 6]{KNU98}. In addition to being simple, these rules are memory-free, i.e., require no memory for their operation.

List organizing algorithms have been analyzed over the past fifty years, e.g., see review on self-organizing linear search in~\cite{HEH85}. While the literature on the move-to-front rule (and the corresponding least-recently-used caching algorithm) is extensive (see, e.g.,~\cite{FIL96b,FIL96,FGT92,JEL97,SLT85a,BP04} and references therein), the results on the transposition rule are scarce. Early analysis of the transposition rule can be found in~\cite{RIV76}. In the same paper it was conjectured that the rule is optimal with respect to the expected value of the search cost. However, it was shown in~\cite{ANW82} that this conjecture is not true in general. Except for the papers mentioned above, the probabilistic analysis of the transposition rule is either limited to the case of simplistic distributions~\cite{KRo80,TNe82} or numerical studies~\cite{BEY97,LLS84}. The reader is referred to~\cite{SLT85a} for a combinatorial (amortized) analysis of the transposition rule.

In the present paper we develop a new approach for analyzing the transposition rule. The approach is based on a coupling with a constrained asymmetric exclusion process. This allows us to establish an asymptotic optimality of the rule for two families of request distributions. Specifically, we prove that the logarithm of the tail probability of the search cost is asymptotically optimal under the transposition rule when the request distribution is either power law or geometric.

The rest of the paper is organized as follows. The model description and main results can be found in the next section. In Section~\ref{sec:AEP} we describe an associated asymmetric exclusion process and characterize its stationary behavior. Section~\ref{sec:tr} relates the exclusion process to the transposition rule for self-organizing lists. Section~\ref{sec:mainresult} contains the proofs of the results stated in Section~\ref{sec:ModelMainResult}. Conclusions and some open questions are discussed in Section~\ref{sec:Conclusions}.

\section{Model and results}\label{sec:ModelMainResult}

We consider an infinite list of items $L=\{1,2,\ldots,N,\ldots\}=\mathbb{N}$. At integer times $t=0,1,2,\ldots$ a request arrives for an item from $L$. The item requested at time $t$ is denoted by $R(t)$. The requests are independent and identically distributed, and $\pi_i$ denotes the probability of item $i$ being requested, $\sum_{i\geq 1}\pi_i=1$. Without loss of generality we assume that $\pi_i \geq \pi_{i+1}$ for all $i$. Let $R$ be equal in distribution to $R(t)$, i.e., $\Pr[R=i] = \pi_i$.

The evolution of the list $L$ is governed by the transposition rule. At time $t=0$ the list is assumed to be ordered as $\{1,2,\ldots,N,\ldots\}$. Upon every request, the requested object is moved forward by one position in the list while the object in front of it is moved one position back. If the first item in $L$ is requested, the list does not change. The basic idea is that frequently requested items are moved closer to the beginning of the list over time; on the other hand, items with low request probabilities end up at some distance from the beginning of the list.

At every time $t$ the list is represented as some permutation $\sigma:\mathbb{N}\rightarrow \mathbb{N}$. Let $X_{i}(t)$ be the position of the item $i$  in the list at time $t$. Our focus is on the behavior of the position $C(t):=X_{R(t)}(t)$ of the requested item, i.e., the search cost, as $t\rightarrow\infty$. We note that if permutation $\sigma$  is fixed, then the distribution of $C(t)$ is determined completely by $\pi:=\{\pi_i\}_{i=1}^\infty$. In this case we use $C^\sigma$ to denote the random position of the selected element $R$. Namely, $C^\sigma$ is simply the (random) search cost required to locate the requested item in a {\em given} list order $\sigma$. Thus, there exist two sources or randomness affecting the search cost: one corresponding to the random arrangement $\sigma$ of the items, and one corresponding to the randomness of the requested item $R$.

Our first lemma is a simple observation stating that for every permutation $\sigma$ the tail asymptotics of $C^\sigma$ dominates the tail asymptotics of $R$.

\begin{lemma}\label{lemma:LowerBound}
For any distribution $\pi$, permutation $\sigma$, and for every
$x\in\mathbb{N}$
\[
\Pr[C^\sigma>x]\geq \Pr[R>x].
\]
\end{lemma}

\begin{proof} The definition of the search cost renders
\begin{align*}
\Pr[C^\sigma>x] &= \sum_{j:\;\sigma(j)>x} \pi_j \nonumber \\
& \geq \sum_{j>x} \pi_j =\Pr[R>x], 
\end{align*}
where the inequality is holds by the monotonicity of elements of $\pi$.
\end{proof}

Thus, as far as the tail probability asymptotics is concerned, no list ordering algorithm can achieve a better performance than the one under the optimal static arrangement. Note that arranging item in the decreasing order of $\pi_i$ is feasible only if the distribution $\pi$ is known in advance.

We say that $R$ is distributed as a power law with parameter $\alpha>1$ if $\pi_i =ci^{-\alpha}$ for all $i$, where $c^{-1}=\sum_{i\geq 1}i^{-\alpha}$ is the normalization constant. Random variable $R$ is defined to be asymptotically geometric with parameter $0<\nu<1$ when $i^{-1} \log\pi_i \to \log\nu$ as $i \to \infty$. The next result states that the transposition rule is asymptotically optimal with respect to the logarithm of the tail asymptotics for these two distribution families.

\begin{theorem}\label{thm:MainResult}
Let $\pi$ be either power law with parameter $\alpha>1$ or asymptotically geometric with parameter $0<\nu<1$. Then
\[
\limsup_{x\rightarrow\infty}\limsup_{t\rightarrow\infty}{\log\Pr[C(t)>x]\over \log\Pr[R>x]}=1. 
\]
\end{theorem}

\begin{proof}
See Section~\ref{sec:mainresult}.
\end{proof}

Often it is of interest to consider list that contain only a finite number of items, i.e., $\pi$ has a finite support. Although we will not make use of the following fact, we remark that for every distribution $\pi$ with finite support on $\{1,2,\ldots,N\}$ ($\pi_i=0$ for all $i>N$) the described system is an irreducible, reversible, aperiodic Markov chain, and the unique stationary solution is of the following product form
\[
\Pr[X_{1}=i_1,\; X_{2}=i_2, \ldots, X_{N}=i_N] = \frac{\prod_{j=1}^N \pi_j^{-i_j}}{\sum_{(k_1,\ldots,k_N) \in P_N} \prod_{j=1}^N \pi_j^{-k_j}},
\]
where $P_N$ denotes the set of all permutations of list $L_N = \{1,2,\ldots,N\}$. A natural way to introduce a power law and geometric distribution for the case of finite support is to take the distribution $\pi$ conditioned on event $\{i\leq N\}$. Denote by $\pi_N$ the truncated distribution and let random variable $R_N$ be defined by $\Pr[R_N>x]:=\Pr[R>x\,|\,R\leq N]$. Note that the existence of a unique stationary distribution for every $N$ allows us to consider the stationary search cost denoted by $C_N$.

\begin{theorem}\label{thm:MainResultTruncated}
Let either (i) $\pi_N$ be truncated power law with parameter $\alpha>1$ and $x/N < \gamma$ for some $\gamma<1$ or (ii) $\pi_N$ be truncated asymptotically geometric with parameter $0<\nu<1$ and $x<N$. Then 
\[
\lim_{\{x,N\}\rightarrow\infty}{\log\Pr[C_N>x]\over \log\Pr[R_N>x]}=1. 
\]
\end{theorem}

\begin{proof} See Section~\ref{sec:mainresult}.
\end{proof}

\section{Constrained asymmetric exclusion process}
\label{sec:AEP}

\begin{figure}[t]
\centering \epsfig{figure=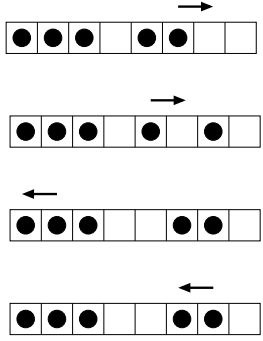, height=1.75in} \caption{An example of evolution of the system with 4 particles. The arrows indicate intended movements of particles. Movements actually occur only in the first two instances.} \label{fig:caep_eg}
\end{figure}

In this section we consider a certain constrained asymmetric exclusion process in which we examine the deviation of the boundary particle from its minimal position. In particular, we consider an $n$ particle system on countably many slots on a half-line enumerated from left to right as $1,2,\ldots$. Each particle is associated with an independent Poisson process of unit intensity (the actual rate is not important). At arrival instances of their corresponding Poisson processes particles move left or right with probabilities $p$ and $q$, respectively. Multiple occupancies are not allowed and a move actually occurs only if the target slot is empty (see Fig.~\ref{fig:caep_eg} for an example). A particle can not move left if it is located in the first slot.

Assume that $p>q$ and define $\beta := q/p < 1$. Given that particles are enumerated from left to right with natural numbers, let $Z_i$ be the position (slot number) of the $i$-th particle.  We first  verify in a straightforward way that the stationary distribution is of the
following form
\begin{equation}
\Pr\left[Z_1=i_1, \, Z_2=i_2, \ldots, Z_n=i_n \right] = \eta_{n}^{-1} \, \beta^{\sum_{j=1}^n i_j}, \label{eq:reversibleCAEP}
\end{equation}
for all  $1 \leq i_1<i_2<\cdots<i_n$, where $\eta_{n}$ is the normalization constant. To this end
\begin{align*}
\sum_{1\leq i_1<i_2<\cdots<i_n}\beta^{\sum_{j=1}^n i_j} & \leq \sum_{1\leq j\leq n}\sum_{1\leq i_j<\infty}\beta^{\sum_{j=1}^n i_j} \\
&= \left(\frac{\beta}{1-\beta}\right)^n<\infty,
\end{align*}
and, therefore, the normalization constant $\eta_n$ is finite. It is easy to check that the underlying Markov chain is irreducible, reversible, aperiodic and that~(\ref{eq:reversibleCAEP}) satisfies the stationary equation. Hence,~(\ref{eq:reversibleCAEP}) indeed describes the stationary distribution.

We point out that the minimal possible value of $\sum_{i=1}^n Z_i$ is $\sum_{i=1}^n i = n(n+1)/2$. Throughout the paper we interpret $\prod_{i=j}^k (\cdot)\equiv 1$ for $k<j$.

\begin{lemma} The normalization constant satisfies
\[
\eta_{n} = \beta^{\frac{n(n+1)}{2}} \prod_{i=1}^n \frac{1}{1-\beta^i}.
\]
\label{lemma:eta}
\end{lemma}

\begin{proof}
For each integer $k \geq 0$ let $\eta_{n,k}$ denote the sum of $\beta^{\sum_{j=1}^n i_j}$ over the feasible choices of $i_j$ such that $\max_{1\leq j\leq n}\{i_j\} \leq n+k$. Then clearly $\eta_{n,0}\leq \eta_{n,1}\leq\cdots$ and $\eta_{n,k} \to \eta_n$ as $k \to \infty$. We claim that for every $n,k$ 
\be{eq:eta_nk}
\eta_{n,k} = \beta^{\frac{n(n+1)}{2}} \prod_{i=1}^k
\frac{1-\beta^{n+i}}{1-\beta^i}. 
\ee 
The expression for $\eta_n$ follows immediately by taking the limit as $k\rightarrow\infty$ in~(\ref{eq:eta_nk}).

The proof of~(\ref{eq:eta_nk}) is by induction. It is trivial to validate that $\eta_{n,0} = \beta^{\frac{n(n+1)}{2}}$, $n \geq 1$,
and $\eta_{1,k} = \sum_{i=1}^{k+1} \beta^i$, $k \geq 0$, conform to~(\ref{eq:eta_nk}). Next, we assume that~(\ref{eq:eta_nk}) holds for $\eta_{i,j}$ for all $i,j$ such that either $i\leq n,\, j<k$ or $i < n,\, j \leq k$ and show that the statement is true for $\eta_{n,k}$. The quantity $\eta_{n,k}$
satisfies the following equality
\begin{equation}
\eta_{n,k} = \eta_{n,k-1} + \beta^{n+k} \eta_{n-1,k};
\label{eq:eta1}
\end{equation}
the first term corresponds to the case $\max_{1\leq j\leq n}\{i_j\} < n+k$, while the second one to the case $\max_{1\leq
j\leq n}\{i_j\}= n+k$, i.e., the last ($n$-th) particle is in the slot $n+k$. From~(\ref{eq:eta1}) and the inductive assumption one derives
\begin{align*}
\eta_{n,k} &= \beta^{\frac{n(n+1)}{2}} \prod_{i=1}^{k-1}
\frac{1-\beta^{n+i}}{1-\beta^i} + \beta^{\frac{(n-1)n}{2}+n+k}
\prod_{i=1}^k \frac{1-\beta^{n-1+i}}{1-\beta^i} \\
&= \left(1+ \beta^k \frac{1-\beta^n}{1-\beta^k} \right)
\beta^{\frac{n(n+1)}{2}} \prod_{i=1}^{k-1}
\frac{1-\beta^{n+i}}{1-\beta^i} \\
&=\beta^{\frac{n(n+1)}{2}} \prod_{i=1}^{k}
\frac{1-\beta^{n+i}}{1-\beta^i}.
\end{align*}
This concludes the proof.
\end{proof}

Next, we use~(\ref{eq:reversibleCAEP}) to examine the stationary deviation $\kappa_{n}$ of the last particle from its minimal position, i.e., $\kappa_{n} := Z_n - n\geq 0$. Expressions~(\ref{eq:reversibleCAEP}),~(\ref{eq:eta_nk}) and Lemma~\ref{lemma:eta} yield
\begin{align}
\Pr[\kappa_{n} = i] &= \frac{\eta_{n-1,i} \, \beta^{n+i}}{\eta_{n}} \nonumber\\
&= \beta^i (1-\beta^n) \prod_{j=1}^{n-1} (1-\beta^{i+j})<\beta^i, \label{eq:kappa_nk}
\end{align}
and, thus, 
\be{eq:geomtailKappan} \Pr[\kappa_{n} \geq i]< (1-\beta)^{-1} \beta^i. 
\ee 
Interestingly, this implies that there exists a limiting behavior for the case when the number of particles $n$ grows to infinity. Indeed, as $n\rightarrow\infty$, the random variable $\kappa_n$ converges in distribution to a random variable $\kappa$ with distribution given by 
$$
\Pr[\kappa=i]=\beta^i\prod_{j=1}^{\infty}(1-\beta^{i+j}). 
$$ 
From the preceding we conclude that $\kappa$ is asymptotically geometric with parameter~$\beta$
\begin{equation*}
\lim_{i \to \infty} \Pr[\kappa=i] \beta^{-i} = 1,
\end{equation*}
and that random variable $\kappa$ is stochastically monotone in parameter $\beta$. Finally, we note that for $\beta<1/2$, or equivalently $2q < p$, the most probable value of $\kappa$ is zero.

\section{Coupling}
\label{sec:tr}

The following lemma relates the stationary (as $t \to \infty$) properties of list $L$ operating under the transposition rule and characteristics of the particle system studied in Section~\ref{sec:AEP}. Let $\kappa_i(\beta)$ explicitly denote the dependency of the random variable $\kappa_i$ on parameter $\beta$ (see~(\ref{eq:kappa_nk})) and set $\beta:= \pi_{i+1}/\pi_i \leq 1$.

\begin{proposition}\label{proposition:connection} For every  $x\geq i \geq 1$
\[
\limsup_{t\to\infty} \; \Pr\left[\bigvee_{j=1}^i X_{j}(t) > x
\right] \leq \Pr\left[\kappa_i(\beta)+i
> x\right].
\]
\label{lemma:connection}
\end{proposition}

\begin{remark}
When the support of $\pi$ is finite there exists a unique stationary solution and, therefore, the left side of the preceding inequality converges as $t \to \infty$.
\end{remark}

\begin{proof}
The proof is based on a coupling argument. We start by exploiting a Poisson embedding technique (see~\cite{FIH96} for an application of the technique in the context of the move-to-front rule). The requests for item $i$, form a Poisson process of intensity $\pi_{i}$. Then, the limiting behaviors (as $t\to\infty$) of the original discrete-time system and the system with the Poisson request patterns are the same.

Given a Poisson process (set of arrival times) $\Lambda$ with rate $\lambda$, let $\Lambda(p)$ denote its subset, $\Lambda(p)
\subseteq \Lambda$, formed by including each element of $\Lambda$ in $\Lambda(p)$ independently with probability $p$. Let $\Lambda_i$ be the set of request times for item $i \in \mathbb{N}$.

Next we construct a modified list $\hat L$ consisting of the same items as the original list~$L$. Parameters of the new system are
denoted with the "hat" symbol. Specifically, $\hat X_j(t)$ denotes the position of element $j$ at time $t$ in the list $\hat L$. Each element $j \in \hat L$ is associated with an independent Poisson process $\hat \Lambda_j$ defined as
\[
\hat \Lambda_j := \begin{cases} \Lambda_j(\pi_{i}/\pi_{j}), & 1 \leq j \leq i,\\
\Lambda_j \cup \Lambda^+_j, & j>i,
\end{cases}
\]
where $\Lambda^+_j$ is an independent Poisson process with rate $\pi_{i+1} - \pi_{j}$. Note that processes $\hat \Lambda_j$ are constructed in such a way that they are Poisson with rates $\pi_{i}$ for $1 \leq j \leq i$ and $\pi_{i+1}$ for $j>i$. In addition, observe that $\hat \Lambda_j \subseteq \Lambda_j$ for $1 \leq j \leq i$ and $\hat \Lambda_j \supseteq \Lambda_j$ for $j>i$. Furthermore, let function $\varphi_j(t)$ be defined as follows
\[
\varphi_j(t) :=
\begin{cases}
\ind_{\{\exists k>i:\; \hat X_k(t)=\hat X_j(t)-1 \}}, & j \leq i, \\
\ind_{\{\exists k \leq i:\; \hat X_k(t) = \hat X_j(t)-1\}}, & j>i,
\end{cases}
\]
i.e., the function $\varphi_j(t)$ indicates whether item $j$ is preceded by an item $k$ such that rates of $\hat \Lambda_j$ and
$\hat \Lambda_k$ differ. The request process $\hat\Lambda$ to list $\hat L$ is a superposition of Poisson processes
$\hat\Lambda_j$:
\begin{equation}
\hat\Lambda(t) := \bigcup_{j:\; \varphi_j(t)=1} \hat\Lambda_j (t).
\label{eq:unionL}
\end{equation}
Item $j$ is requested from $\hat L$ at time $t=T$ if $T \in \hat\Lambda_j \cap \hat \Lambda$. In other words, item $j\leq i$ is requested according to $\hat\Lambda_j$ only if it is preceded by an item in $\{i+1, i+2, \ldots\}$. On the other hand, requests for item $j>i$ are placed according to $\hat\Lambda_j$ only if element $j$ is preceded by an item in $\{1,\ldots,i\}$. Note that the set of Poisson processes included in $\hat \Lambda$ changes with the evolution of list $\hat L$. In addition, the number of elements in the union in~(\ref{eq:unionL}) is always finite and bounded from above by $2i$.

The modified list $\hat L$ operates under the transposition rule with one modification. Namely, after the transposition
rule rearranges items in either of the lists ($L$ or $\hat L$), a reordering operator ${\cal R}_i$ is applied to $\hat L$. The operator works as follows. The list is divided in two sublist: $\{1,2,\ldots,i\}$ and $\{i+1,i+2,\ldots\}$. The operator ${\cal R}_i$ reorders each sublist of $\hat L$ so that the order of elements within the sublists is the same as in the original list $L$. However, only items belonging to the same sublist are allowed to exchange positions in the list. An example of how ${\cal R}_i$ operates is shown in Fig.~\ref{fig:order_eg}.

\begin{figure}[t]
\centering \epsfig{figure=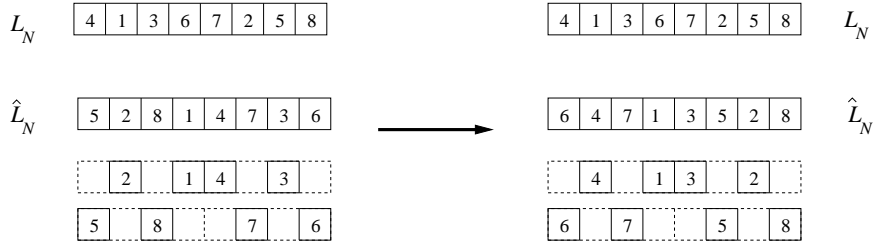, height=1.25in} 
\caption{An example of reordering by operator ${\cal R}_i$ (in this case $i=4$ and $N=8$). The initial states of the lists are shown on the left. The modified list $\hat L_N$ is divided into two sublists: $\{1,2,3,4\}$ and $\{5,6,7,8\}$. Upon reordering items in each of the sublists according to $L_N$, the new ordering in $\hat L_N$ is shown on the right.} \label{fig:order_eg}
\end{figure}

Next, assume that both lists are in the same permutation at time $t=0$. Then we argue that

\begin{lemma}\label{lemma:coupling}
For every $t\in\mathbb{R}_+$ and $1 \leq j \leq i$
\begin{equation}
X_{j} (t) \leq \hat X_{j}(t). \label{eq:X(t)hatX(t)}
\end{equation}
\end{lemma}

\begin{proof}
Either list changes only at times of requests to one of the systems; denote those times by $0<T_1<T_2<\cdots<T_n<\cdots$. Since there are no changes in item order between times $\{T_n\}$, it is sufficient to prove that~(\ref{eq:X(t)hatX(t)}) holds for $t= T_n+$, $n \geq 1$. To this end, suppose that~(\ref{eq:X(t)hatX(t)}) holds for $t=T_{n-1}+$ and consider the two lists at time $t=T_{n}+$. There are three cases that need to be examined:

\begin{itemize}
\item[(i)] {\em At time $t=T_{n}$ item $j \leq i$ is requested from $L$.} By the construction of $\hat \Lambda$, this event implies
that item $j$ is requested in the modified list $\hat L$ only with some probability depending on the state of $\hat L$.

If item $j$ is not requested from $\hat L$, then the set of positions occupied by items $\{1,\ldots,i\}$ in $\hat L$ remains the same. On the other hand, the set of positions occupied by items $\{1,\ldots,i\}$ in $L$ either does not change (when an item in $\{1,\ldots,i\}$ precedes item~$j$ or $j$ is the first item in the list) or is pairwise smaller (when an item in $\{i+1,\ldots,N\}$ precedes item~$j$). This, in conjunction with the fact that the relative order of items $1,\ldots,i$ in $L$ and $\hat L$ at $t=T_n+$ is the same (due to ${\cal R}_i$), implies that~(\ref{eq:X(t)hatX(t)}) holds for $t=T_n+$.

On the other hand, if item $j$ is requested in both lists, then the only case that needs to be examined in detail is the one in which $j$ is preceded by item $l \leq i$ in $L$; note that from~(\ref{eq:unionL}) it follows that $j$ is preceded by an item $v>i$ in $\hat L$. Since $X_l(T_n-) \leq \hat X_l(T_n-)$ and the order of the items $1,\ldots,i$ is the same in both lists (due to ${\cal R}_i$) one has 
\be{eq:jl} 
X_j(T_n-)+1 \leq \hat X_j(T_n-),
\ee
The fact that items $j$ and $l$ are transposed in $L$ but not in $\hat L$ yields that the order of $j$ an $l$ is different in the two lists before
${\cal R}_i$ is applied. Thus, ${\cal R}_i$ exchanges positions of items $j$ and $l$ in~$\hat L$. This leads to
\[
\hat X_j(T_n+) = \hat X_l (T_n-) \geq X_l(T_n-) = X_j(T_n+)
\]
and
\[
\hat X_l(T_n+) = \hat X_j(T_n-) - 1 \geq X_j (T_n-) = X_l(T_n+),
\]
where the first inequality follows from the inductive assumption and the second inequality is due to~(\ref{eq:jl}). Thus, we conclude that~(\ref{eq:X(t)hatX(t)}) holds for $t=T_n+$.

\item[(ii)] {\em At time $t=T_{n}$ item $j>i$ is requested from $\hat L$.} The argument is very similar to the one in~(i). In this case $j$ is preceded by some $v\leq i$ in $\hat L$ (see~(\ref{eq:unionL})).

If item $j$ is not requested from $L$ or is preceded (in $L$) by an item $v>i$, then the positions occupied by items $1,\ldots,i$ in $L$ do not change. However, in $\hat L$ item $j$ must be preceded by $v\leq i$ and, thus, $v$ is moved one position back by the transposition rule.
Alternatively, if in $L$ item $j$ is requested and preceded by $l \leq i$ then either $v=l$ and $X_l(T_n-) = \hat X_v(T_n-)$ or $v \not= l$ and 
\[
X_l(T_n-) \leq \hat X_l(T_n-) + 1,
\] 
\[
X_v(T_n-) \leq \hat X_v(T_n-) + 1.
\]
In either case after the items are transposed one has $X_l(T_n+) \leq \hat X_l(T_n+)$ and $X_v(T_n+) \leq \hat X_v(T_n+)$.

\item[(iii)] {\em At time $t=T_n$ item $j>i$ requested from $L$ but not from $\hat L$.} This implies that $j$ is preceded in $\hat L$ by an
item $v>i$ (see~(\ref{eq:unionL})). If $j$ is preceded in $L$ by an item $l>i$ the positions occupied by items $1,\ldots,i$ do not change in either list. Therefore, we only need to consider the case when $j$ is preceded by $l \leq i$. However, in that case we necessarily have
\[
\hat X_l(T_n-) \geq X_l(T_n-) +1,
\]
which implies
\[
X_l(T_n+) = X_l(T_n-) +1 \leq \hat X_l(T_n-) = \hat X_l(T_n+).
\]
\end{itemize}
The preceding (i)-(iii) establish~(\ref{eq:X(t)hatX(t)}).
\end{proof}

Next, let variables $Z_j(t)$, $1\leq j \leq i$, be defined by
\[
(Z_1(t),\ldots, Z_i(t)) := {\cal S}(\hat X_1(t), \ldots, \hat
X_i(t)),
\]
where $\cal S$ is the sorting (in the increasing order) operator. Note that the definition of $Z_i(t)$ and~(\ref{eq:X(t)hatX(t)}) imply
\begin{equation}
Z_i(t) \geq \bigvee_{j=1}^i  X_j(t). \label{eq:ZX}
\end{equation}
Observe that the evolution of $\{Z_j(t)\}$ is probabilistically the same as in the constrained particle system described in Section~\ref{sec:AEP} (recall that $Z_j(t)$ denotes the position of the $j$th particle at time $t$ in Section~\ref{sec:AEP}) with $p=\pi_i/(\pi_i+\pi_{i+1})$ and $q=\pi_{i+1}/(\pi_i+\pi_{i+1})$. Indeed, $Z_j(t)$ increases by one at Poisson rate $\pi_{i+1}$ only if $Z_{j+1}(t) \not= Z_j(t)+1$, and it decreases by one at rate $\pi_i$ only if $Z_{j-1}\not=Z_j(t)-1$. 

Taking maximums on both sides of~(\ref{eq:ZX}) and applying operator $\Pr~[~\cdot~>~x~]$ leads to
\begin{align*}
\Pr\left[\bigvee_{j=1}^i X_{j}(t) > x \right] &\leq \Pr\left[Z_{i}(t) > x \right].
\end{align*}
From the preceding inequality one obtains
\begin{align*}
\limsup_{t\rightarrow\infty}\Pr\left[\bigvee_{j=1}^i X_{j}(t) >
x\right] &\leq \Pr\left[Z_i
> x \right] \\
&= \Pr[k_i(\beta) + i > x],
\end{align*}
where the last equality follows from the definition of $\kappa_i(\beta)$ in Section~\ref{sec:AEP} and $\beta = \pi_{i+1}/\pi_i$.
\end{proof}

\section{Proofs}\label{sec:mainresult}

Proposition~\ref{proposition:connection} is the primary tool in establishing our results on the performance of the transposition rule. The following lemma is a simple consequence of Proposition~\ref{proposition:connection}.

\begin{lemma} For any $y \geq 1$ and distribution of requests $\pi$
\[
\limsup_{t \to \infty}\Pr[C(t)>x] \leq \Pr[\kappa_y(\pi_{y+1}/\pi_{y})>x-y] + \Pr[R>y,\, R+\kappa_R(\pi_{R+1}/\pi_{R})>x].
\]
\label{lemma:upper}
\end{lemma}

\begin{proof} Conditioning on the requested item and using the monotonicity of the $\max$-operator result in
\begin{align*}
\limsup_{t\to\infty}\Pr\left[C(t)>x \right] &= \limsup_{t\to\infty}\left(\sum_{i\geq 1} \pi_{i} \, \Pr[X_{i}(t)>x]\right) \\
&\leq \sum_{i=1}^y \pi_{i}
\limsup_{t\to\infty}\Pr\left[\bigvee_{j=1}^{y} X_{j}(t)> x\right]
+ \sum_{i\geq y+1} \pi_{i}  \limsup_{t\to\infty}
\Pr\left[\bigvee_{j=1}^{i} X_{j}(t) > x \right] \\
&\leq \Pr[\kappa_y(\pi_{y+1}/\pi_{y})
>x-y] + \Pr[R>y,\, R+\kappa_R(\pi_{R+1}/\pi_{R})>x],
\end{align*}
where the last inequality follows from Proposition~\ref{proposition:connection}.
\end{proof}

At this point we present the proofs of Theorems~\ref{thm:MainResult} and~\ref{thm:MainResultTruncated}.

\begin{proof}[Proof of Theorem~\ref{thm:MainResult}]
The lower bound is an immediate consequence of Lemma~\ref{lemma:LowerBound} and holds for any distribution of requests $\pi$. Hence, we only consider the upper bound.

We first examine the case when $\pi$ is asymptotically geometric with parameter $\nu$. Fix arbitrary small $\epsilon>0$ such that $\nu+\epsilon<1$. By the assumption, there exists $i_\epsilon$ such that $\nu-\epsilon<\pi_{i+1}/\pi_{i}<\nu+\epsilon$ for all $i \geq i_\epsilon$. For any $s^{-1}> \nu+\epsilon$, setting $y=i_\epsilon$ in Lemma~\ref{lemma:upper} yields
\begin{align}
\limsup_{t\to\infty}\Pr[C(t)>x] &\leq \Pr[\kappa_{i_\epsilon}(\nu+\epsilon)>x-i_\epsilon] + \Pr[R+\kappa_R(\nu+\epsilon) >x] \nonumber \\
&\leq\Pr[s^{\kappa_{i_\epsilon}(\nu+\epsilon)}>s^{x-i_\epsilon}] + \Pr[s^{R+\kappa_R(\nu+\epsilon)} > s^x] \nonumber \\
&\leq s^{-x} s^{i_\epsilon} \expect s^{\kappa_{i_\epsilon}(\nu+\epsilon)} +s^{-x} \, \expect s^{R + \kappa_R(\nu+\epsilon)}, \label{eq:geoin2}
\end{align}
where the last step is due to Markov's inequality. From $s^{-1}>\nu+\epsilon$ and (\ref{eq:geomtailKappan}) it follows
$\expect s^{R + \kappa_R(\nu+\epsilon)} < \infty$. This bound,~(\ref{eq:geomtailKappan}) and~(\ref{eq:geoin2}) result in
\begin{align}
x^{-1} \limsup_{t\to\infty}\log \Pr[C(t)>x] &\leq - \log s + x^{-1} \log \left(s^{i_\epsilon} \expect s^{\kappa_{i_\epsilon}(\nu+\epsilon)} +
\expect s^{R + \kappa_{R}(\nu+\epsilon)}\right) \nonumber \\
&\to - \log s, \label{eq:geombound1a}
\end{align}
as $x \to \infty$. On the other hand, note that
\begin{align*}
\Pr[R>x] &\geq
\sum_{i>x}\pi_{i_\epsilon}(\nu-\epsilon)^{i-i_\epsilon} \\
&\geq \pi_{i_\epsilon}(1-\nu+\epsilon)^{-1}(\nu-\epsilon)^{x+1-i_\epsilon},
\end{align*}
implying
\begin{equation}
\liminf_{x\to \infty} x^{-1}\log\Pr[R>x]\geq \log(\nu-\epsilon).
\label{eq:geombound2a}
\end{equation}
Combining~(\ref{eq:geombound1a}) and~(\ref{eq:geombound2a}), we obtain
\[
\limsup_{x \to \infty} \limsup_{t\to\infty} \frac{\log
\Pr[C(t)>x]}{\log\Pr[R>x]}\leq \frac{\log
s^{-1}}{\log(\nu-\epsilon)};
\]
passing $s^{-1} \downarrow \nu + \epsilon$ and then letting $\epsilon \downarrow 0$ yield the result.

Next, we consider the case when $\pi$ is power law with parameter $\alpha>1$. Lemma~\ref{lemma:upper} and~(\ref{eq:geomtailKappan}) yield
\begin{align*}
\limsup_{t\to\infty}\Pr[C(t)>x] &\leq \Pr[\kappa_y(\pi_{y+1}/\pi_y) > x-y] + \Pr[R>y] \\
&\leq \left(1-\frac{\pi_{y+1}}{\pi_{y}} \right)^{-1} \left(\frac{\pi_{y+1}}{\pi_{y}}\right)^{x-y} + \Pr[R>y].
\end{align*}
Letting $y=\lceil \epsilon x/\log x\rceil$ for a sufficiently small $\epsilon>0$ results in an estimate on the two terms in the preceding sum, as $x \to \infty$,
\[
\left(1-\frac{\pi_{y+1}}{\pi_{y}} \right)^{-1} \left(\frac{\pi_{y+1}}{\pi_{y}}\right)^{x-y} = \alpha^{-1} e^\alpha y x^{-\alpha/\e} (1+o(1))
\]
and
\[
\Pr[R>y]=\sum_{i>y}ci^{-\alpha}=O(y^{-\alpha+1}).
\]
Therefore, as $x\to\infty$,
\[
\limsup_{t\to \infty}\Pr[C(t)>x] \leq \Pr[R>y] (1+o(1)).
\]
The preceding equation together with the fact that $\pi$ is a power law yield the statement of the statement of the theorem.
\end{proof}

\begin{proof}[Proof of Theorem~\ref{thm:MainResultTruncated}]
The proof is very similar to the one of Theorem~\ref{thm:MainResult} and, thus, we omit details. Since
\[
\frac{\log \Pr[C_N>x]}{\log \Pr[R_N>x]} = \frac{\log \Pr[C_N>x]}{\log \Pr[R>x]} \frac{\log\Pr[R>x]}{\log\Pr[R_N>x]}
\]
and the upper bound on $\Pr[C_N>x]$ is the same as the one on $\limsup \Pr[C(t)>x]$, we only need to verify that the last fraction in the preceding equality tends to one as $\{x,N\} \to \infty$. However, that easily follows from the assumptions of the theorem.
\end{proof}

\section{Concluding remarks}\label{sec:Conclusions}

We presented an analysis of the transposition rule based on a coupling with a constrained exclusion process. As an outcome, we established an asymptotic optimality of the transposition rule in linear lists. Specifically, when the probability distribution of the requests is power law or geometric we showed that, under the transposition rule, the logarithm of the tail probability of the search cost is asymptotically optimal. 

While the steady-state distribution of the search cost is a primary quantity of interest, rates of convergence play an important role in assessing the applicability of self-organizing algorithms in practice. The proposed coupling may offer new directions for understanding these rates under the transposition rule. The same question, for the related move-to-front algorithm, was investigated in~\cite{FIL96}. As pointed out in~\cite{FIL96}, the transposition rule is expected to have slower rates of convergence than the move-to-front rule. 

%


\small

\begin{thebibliography}{10}

\bibitem{ANW82}
E.~Anderson, P.~Nash, and R.~Weber.
\newblock A counterexample to a conjecture in optimal list ordering.
\newblock {\em J. Appl. Probab.}, 19(9):730--732, 1982.

\bibitem{BEY97}
R.~Bachrach and R.~El-Yaniv.
\newblock Online list accessing algorithms and their applications: {R}ecent
  empirical evidence.
\newblock In {\em Proc. ACM-SIAM SODA}, 1997.


\bibitem{BP04}
J.~Barrera and C.~Paroissin.
\newblock On the distribution of the search cost for the move-to-front rule
  with random weights.
\newblock {\em J. Appl. Probab.}, 41(1):250--262, 2004.

\bibitem{FIL96b}
J.~Fill.
\newblock An exact formula for the move-to-front rule for self-organizing
  lists.
\newblock {\em J. Theoret. Probab.}, 9(1):113--159, 1996.

\bibitem{FIL96}
J.~Fill.
\newblock Limits and rate of convergence for the distribution of search cost
  under the move-to-front rule.
\newblock {\em Theoret. Comput. Sci.}, 164:185--206, 1996.

\bibitem{FIH96}
J.~Fill and L.~Holst.
\newblock On the distribution of search cost for the move-to-front rule.
\newblock {\em Random Structures Algorithms}, 8(3):179-186, 1996.

\bibitem{FGT92}
P.~Flajolet, D.~Gardy, and L.~Thimonier.
\newblock Birthday paradox, coupon collector, caching algorithms and
  self-organizing search.
\newblock {\em Discrete Appl. Math.}, 39:207--229, 1992.

\bibitem{HEH85}
J.~Hester and D.~Hirchberg.
\newblock Self-organizing linear search.
\newblock {\em Computing Surveys}, 17(3):295--311, September 1985.

\bibitem{JEL97}
P.~Jelenkovi\'{c}.
\newblock Asymptotic approximation of the move-to-front search cost
  distribution and least-recently-used caching fault probabilities.
\newblock {\em Ann. Appl. Probab.}, 9(2):430--464, 1999.

\bibitem{KRo80}
Y.~Kan and S.~Ross.
\newblock Optimal list order under partial memory constraints.
\newblock {\em J. Appl. Probab.}, 17:1004--1015, 1980.

\bibitem{KNU98}
D.~Knuth.
\newblock {\em The Art of Computer Programming}, volume 3. Sorting and
  Searching.
\newblock Addison-Wesley, 2nd edition, 1998.

\bibitem{LLS84}
K.~Lam, M.~Y. Leung, and M.~K. Siu.
\newblock Self-organizing files with dependent accesses.
\newblock {\em J. Appl. Probab.}, 21:343--359, 1984.

\bibitem{RIV76}
R.~Rivest.
\newblock On self-organizing sequential search heuristics.
\newblock {\em Comm. ACM}, 19(2):63--67, 1976.

\bibitem{SLT85a}
D.~Sleator and R.~Tarjan.
\newblock Amortized efficiency of list update and paging rules.
\newblock {\em Comm. ACM}, 28(2):202--208, 1985.

\bibitem{TNe82}
A.~Tenenbaum and R.~Nemes.
\newblock Two spectra of self-organizing sequential search algorithms.
\newblock {\em SIAM J. Comput.}, 11(3):557--566, 1982.

\end{thebibliography}

\end{document}